\documentclass[12pt,twoside,leqno]{article}
\usepackage[english]{babel}
\usepackage{amssymb, amsthm, amsmath}
\usepackage{graphicx}

\newcommand{\tr}{\mathrm{tr}} 
\newcommand{\re}{\mathop{\mathrm{Re}}} 
\newcommand{\ii}{\mathrm{i}} 
\newcommand{\Lf}{\mathrm{L}} 
\newcommand{\dd}{\;\mathrm{d}}
\newcommand{\e}{\mathrm{e}}
\newcommand{\pf}{\noindent {\bf PROOF.} \quad}

\newcommand{\C}{\mathbb{C}}

\newcommand{\Z}{\mathbb{Z}}
\newcommand{\PP}{\mathbb{P}}

\newcommand{\TT}{\mathbb{T}}

\newtheorem{thm}{Theorem}
\newtheorem{defn}[thm]{Definition}

\newtheorem{lem}[thm]{Lemma}

 \newcommand{\addressa}{Oliver T. Dasbach 
 \\
 Department of Mathematics
 \\
 Louisiana State University
 \\
 Baton Rouge, LA 70803, USA
 }
 \newcommand{\emaila}{\texttt{kasten@math.lsu.edu}}

\newcommand{\addressb}{Matilde N. Lal\'{\i}n
\\
Department of Mathematical and Statistical Sciences
\\
University of Alberta
\\
Edmonton, AB T6G 2G1, Canada
}
\newcommand{\emailb}{\texttt{mlalin@math.ualberta.ca}}

 \newcommand{\firstpageR}{1}
\begin{document}
\title{On the recurrence of coefficients in the L\"uck-Fuglede-Kadison determinant}
\author{Oliver T. Dasbach and Matilde N. Lalin}

\maketitle
 \begin{abstract}
In this note, we survey results concerning variations of the L\"uck-Fuglede-Kadison determinant with respect to the base group. 
Further, we discuss recurrences of coefficients in the determinant for certain distinguished base groups.
The note is based on a talk that the second author gave at the ``Segundas Jornadas de Teor\'{\i}a de N\'{u}meros'', Madrid, 2007.
 \end{abstract}
\medskip

The object that we consider in this note is given by the following
\begin{defn} {\bf \cite{Dasbach-Lalin}}
Let $\Gamma$ be a group finitely generated by $x_1, \dots, x_l$. 
Let $P =\sum_{g \in \Gamma} c_g g \in \C \Gamma$ such that $c_g = \overline{c_{g^{-1}}}$. 
Let $\lambda$ be a small complex number. More precisely, $|\lambda| < \frac{1}{l_1(P)}$, the reciprocal of the sum of the absolute values of the coefficients of $P$. 
The Mahler measure (or L\"uck-Fuglede-Kadison determinant \cite{Lueck}) of $ 1 - \lambda P$ is given by
\[ m_\Gamma(P,\lambda) = - \sum_{n=1}^\infty \frac{a_n \lambda^n}{n},\]
where $a_n = [P^n]_0$ is the constant coefficient of  the $n$-th power of $P$; in other words, $a_n$ is the trace of the element $P^n \in \C\Gamma$. 
\end{defn} 


We will often consider the generating function for the $a_n$'s
\[u_\Gamma(P,\lambda) = \sum_{n=0}^\infty a_n \lambda^n.\]

Thus,
\[u_\Gamma(P,\lambda)=a_0-\lambda \frac {\mbox{d}}{\mbox{d}\lambda} m_\Gamma(P,\lambda).\]

\section{Summary of previous results}

We have studied in \cite{Dasbach-Lalin} some properties of $m_\Gamma (P, \lambda)$ where we emphasize the variation of $\Gamma$. 
For particular cases of $\Gamma$ one can find formulas for the Mahler measure over $\Gamma$. In the following formulas, $|\lambda| < \frac{1}{l_1(P)}$.

\begin{itemize}
\item If $\Gamma = \Z^h$, 
\begin{equation}
\re(m_\Gamma(P, \lambda)) = m(1 - \lambda P), 
\end{equation}
where the term on the right hand-side indicates the Mahler measure in the classical sense,
\[ m(1 -\lambda P) = \frac{1}{(2 \pi \ii)^h} \int_{\TT^h} \log|1 -\lambda P(x_1, \dots, x_h)| \frac{\dd x_1}{x_1} \dots \frac{\dd x_h}{x_h}.\]
Here $\TT^h = \{ |x_1|=\dots = |x_h|=1\}$ is the $h$-th unit torus.
\item If $\Gamma$ is finite, 
\begin{equation} \label{coro9}  
m_\Gamma(P, \lambda) = \frac{1}{|\Gamma|} \log \det (1 -\lambda A),
\end{equation}
where $A$ is the adjacency matrix of a weighted Cayley graph of $\Gamma$ generated by the monomials of $P$ (see \cite{Dasbach-Lalin}) and we are considering the main branch of the logarithm. This formula provides a meromorphic continuation of the Mahler measure $m_\Gamma(P, \lambda)$ to the complex plane minus $\mbox{Spec}(A)$.

\item If $\Gamma = \Z/m_1 \Z \times \dots \times \Z/m_h \Z$, then
\begin{equation} \label{finite:eq}
m_\Gamma(P, \lambda) = \frac{1}{|\Gamma|} \log \left(\prod_{j_1, \dots, j_h} (1 -\lambda P(\xi_{m_1}^{j_1}, \dots,\xi_{m_h}^{j_h}))  \right),
\end{equation}
where $\xi_m$ denotes an $m$-th primitive root of the unity, and again, we are considering the main branch of the logarithm. 
\end{itemize}

Of further interest are approximation results for the Mahler measure over infinite groups.
In \cite{Dasbach-Lalin} it is shown that

\begin{itemize}
\item 
\begin{equation}\label{lim}
\lim_{m_i \rightarrow \infty} m_{ \Z/m_1 \Z \times \dots \times \Z/m_h \Z} (P, \lambda) = m_{\Z^h}(P, \lambda)
\end{equation}
\item
\begin{equation}
\lim_{m \rightarrow \infty}m_{D_m}(P, \lambda) = m_{D_\infty}(P, \lambda) 
\end{equation}
Here $D_m$ is the dihedral group  $D_m=\langle \rho, \sigma | \rho^m, \sigma^2, \sigma \rho \sigma \rho \rangle$ and
$D_{\infty}=\langle \rho, \sigma | \sigma^2, \sigma \rho \sigma \rho \rangle$.
\end{itemize}

\section{Recurrence relations}

Following the ideas in \cite{RV:ModularMahler}, notice that if $\Gamma = \Z^l$, then
\[u(\lambda) = \frac{1}{(2\pi\ii)^l} \int_{\TT^l} \omega(\lambda) \qquad \mathrm{where}\quad \omega(\lambda) = \frac{1}{1-\lambda P(x_1, \dots, x_l)} \frac{\dd x_1}{x_1} \dots \frac{\dd x_l }{x_l}.\]
For each value of $\lambda$,  $u(\lambda)$ is a period of $\omega(\lambda)$ in $\PP^{l}(\C)$. The integral depends on the homology class of $\TT^l$ in $H_l(\PP^l(\C) \setminus V, \Z)$, where $V$ is the zero locus of the denominator in the rational function (which is generically non-singular as $\lambda$ varies). See 
Griffiths \cite{Griffiths:Periods}.

Now if we take successive derivatives of $\omega(\lambda)$, we obtain several differential forms that belong to a subspace of the de Rham cohomology $H^l(V)$ which has finite dimension. Griffiths proves that $u(\lambda)$ satisfies a Picard-Fuchs differential equation
\[p_{k}(\lambda) u^{(k)} + p_{k-1}(\lambda) u^{(k-1)}+ \dots + p_0(\lambda) u = 0,\]
where the $p_i$ are polynomials in $\lambda$, see \cite{Griffiths:Periods} for details. 

From such a differential equation it is easy to deduce a linear recurrence with polynomial coefficients for the coefficients $a_n = [P^n]_0$ of $u(\lambda)$.


One can extend this more generally:
\begin{thm}
If $\Gamma$ is a finitely generated abelian group then the coefficients $a_n=[P^n]_0$ satisfy a linear recurrence relation with polynomial coefficients.
\end{thm}
\pf
Let $\Gamma = \Z^l \times \Z/m_1\Z \times \dots \times \Z/m_h \Z$. Combining Equation (\ref{finite:eq}) and the techniques that are used in the proof of Equation (\ref{lim}) (see \cite{Dasbach-Lalin}), we obtain that
\[m_\Gamma(P,\lambda) = \frac{1}{m_1 \dots m_h} \sum_{j_1, \dots, j_h} m(P(x_1,\dots, x_l,\xi_{m_1}^{j_1}, \dots, \xi_{m_h}^{j_h}), \lambda),\]
where $\xi_k$ is a primitive root of unity, and the sum on the right involves Mahler measures in the classical (abelian) sense.

Then the result follows easily since it is known for Mahler measures.

\qed

For finite groups, we have the following
\begin{thm}
If $\Gamma$ is a finite group, then the coefficients $a_n=[P^n]_0$ satisfy a recurrence relation (with constant coefficients) of length at most $|\Gamma|$. 
\end{thm}
\pf In the proof of Equation (\ref{finite:eq}) (Theorem 6 in \cite{Dasbach-Lalin}) we write 
\[a_n =\frac{1}{|\Gamma|} \tr(A^n).\]
Any polynomial that annihilates $A$ yields a recurrence relation with constant coefficients for $a_n$. In particular, the characteristic polynomial yields a bound for the length of the recurrence. 

\qed 

More is known in the case where $\Gamma$ is free \cite{Garoufalidis-Belissard}: the function $u(\lambda)$ turns out to be algebraic. A proof
for this uses algebraic functions in non-commuting variables and a theorem of
Haiman \cite{Haiman}; see \cite{Garoufalidis-Belissard} for details.

A natural question is the following: can we say anything for ``intermediate'' groups?

\section{Some examples}

Of particular interest is the case when $P_l=x_1+x_1^{-1} + \dots + x_l + x_l^{-1}$ where $x_1, \dots, x_l$ are the generators in a given group presentation of the group $\Gamma$.

It is easy to see that $[P_l^n]_0$ has the following interpretation

\begin{lem} \label{circuits}
The number of closed circuits based at the origin in the Cayley graph of $\Gamma$ - with respect to the generators  $x_1, \dots, x_l$ in the presentation - is given by
$[P_l^n]_0$.
\end{lem}






\subsection{Abelian groups}
Consider the example of $P_2=x+x^{-1}+y+y^{-1}$ with group $\Gamma=\Z^2$. Then one obtains (see for example \cite{Domb,RV:ModularMahler})
\[ u(\lambda) = \sum_{n=0}^\infty \binom{2n}{n}^2 \lambda^{2n} .\]

Since we only have terms with even degree, we reparametrize,
\[ v (\mu) = \sum_{n=0}^\infty \binom{2n}{n}^2 \mu^{n}.\] 
Then the differential equation is given by \cite{RV:ModularMahler}
\[\mu(16\mu-1) v'' + (32 \mu -1) v' +4 v =0.\] 

From such a differential equation it is possible to deduce a recurrence of the coefficients, in this case, $a_0=1, a_2=4$, and
\[n^2 a_{2n} -4(2n-1)^2a_{2n-2}=0.\]
Of course this equation can be easily deduced from the formula for $a_n$.

If one considers one more variable, $P_3= x+x^{-1}+y+y^{-1}+z+z^{-1}$. Then one obtains
\[ u(\lambda) = \sum_{n=0}^\infty \binom{2n}{n} \sum_{k=0}^n \binom{n}{k}^2\binom{2k}{k} \lambda^{2n},\] 
and $a_0=1, a_2=6, a_4=90$,
\begin{equation}\label{rec:P_l}
n^3 a_{2n}-2(2n-1)(10n^2-10n+3)a_{2n-2}+36(2n-1)(n-1)(2n-3)a_{2n-4}=0.
\end{equation}

As observed in \cite{Dasbach-Lalin}, for $l$ variables, we obtain
\begin{equation} \label{p1}
a_{2n}^{(P_l)} = \sum_{j_1+\dots+j_l=n} \frac{(2n)!}{(j_1!)^2 \dots (j_l!)^2}.
\end{equation}
By Lemma \ref{circuits}  these coefficients can be interpreted as  the number of circuits of length $2n$ (that start and end at the origin) in the $l$-dimensional cubic lattice.

Another polynomial  that is interesting to study  is 
\[Q_l=(1+x_1+\dots +x_{l-1})(1+x_1^{-1}+ \dots +x_{l-1}^{-1}).\]
One has \cite{Dasbach-Lalin} that
\begin{equation} \label{p2}
a_n^{(Q_l)}=\sum_{j_1+\dots+j_l=n} \left(\frac{n!}{j_1! \dots j_l!}\right)^2.
\end{equation}
We see that the terms in equation (\ref{p1}) correspond to the one in (\ref{p2}) multiplied by $\binom{2n}{n}$, and that allows an easy translation for the recurrences. If
\[ \sum_{k=0}^{l-1} p_k(n) a_{2n-2k}^{(P_l)} =0,\]
then
\begin{eqnarray}
\sum_{k=0}^{l-1} \frac{(n!)^2 (2n-2k)!}{((n-k)!)^2 (2n-2l+2)!} p_k(n) a_{n-k}^{(Q_l)} =0.   \label{recurrsive_relation}
\end{eqnarray}

As an example,
for  $l=3$ we have \[Q_3 = (1+x+y)(1+x^{-1}+y^{-1})=3+ x+x^{-1}+y+y^{-1}+ x y^{-1}+ x^{-1} y.\] 

 One obtains $a_0=1, a_1=15$,
\[n^2a_n -(10n^2-10n+3)a_{n-1}+9(n-1)^2a_{n-2}=0,\]
from Equations (\ref{rec:P_l}) and (\ref{recurrsive_relation}).

Furthermore, from the above discussion a closed form is given by (see also \cite{CFRS,Domb})
\[a_n=\sum_{k=0}^n \binom{n}{k} ^2 \binom{2k}{k},\]
and with the notations of this section:
\[[P_3^{2n}]_0=\binom{2n}{n} a_n = \binom{2n}{n} [Q_3^n]_0.\]

Also note that $b_n:=[P_3^n]_0=[(Q_3-3)^n]_0$ is related to $a_n$ by 
\[b_n=\sum_{j=0}^n   \binom{n}{j}  (-3)^{n-j} a_j.\]

One can show that $b_n$ satisfies the recursion:
$b_0 = 1, b_1 = 0, b_2=6$,
\[n^2 b_n-n(n-1)b_{n-1}-24(n-1)^2b_{n-2}-36(n-2)(n-1) b_{n-3}=0.\]

\begin{figure} 
\includegraphics[width=\textwidth,keepaspectratio]{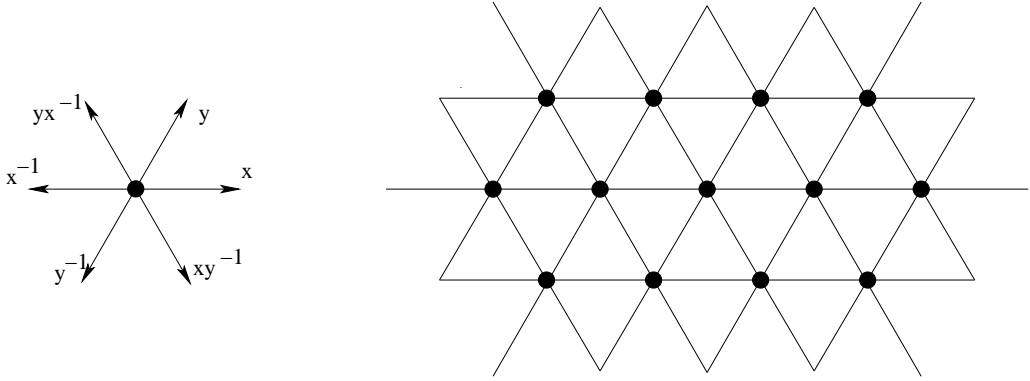}
\caption{The triangular lattice: An edge in the lattice corresponds to a multiplication in $\Z^2$.}
\label{triangular_lattice}
\end{figure}

The numbers $a_n$ and $b_n$ have an interesting interpretation.
Since \[Q_3-3=x+x^{-1}+y+y^{-1}+ x y^{-1}+ x^{-1} y\] the constant coefficient $b_n$ in $(Q_3-3)^n$ counts the number of closed circuits of length $n$, based at a fixed point,  in the triangular lattice that is depicted in Figure \ref{triangular_lattice}.

\begin{figure}
\includegraphics[width=\textwidth,keepaspectratio]{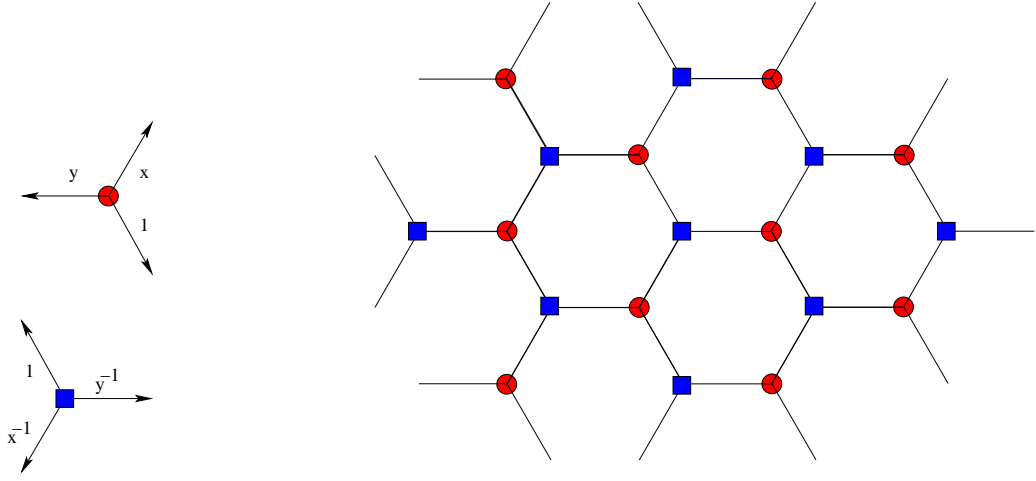}
\caption{The honeycomb lattice.}
\label{honeycomb}
\end{figure}

A look at the honeycomb lattice in Figure \ref{honeycomb} reveals that at any point there are nine different paths of length two originating from that point.
Three are closed paths the other six can be labeled by $x, x^{-1}, y, y^{-1}, x y^{-1}$ and $x^{-1} y$  as in the triangular lattice (Figure \ref{triangular_lattice}).
It follows that $a_n$ is the number of closed circuits in the honeycomb lattice of length $2 n$ that are based at a fixed point.

Our motivation for studying the above examples comes from the fact that
\[m(1+x+y) = \frac 1 2 m((1+x+y)(1+x^{-1}+y^{-1}))\] and
\[m(1+x+y) = \frac{3 \sqrt{3}}{4 \pi} \Lf(\chi_{-3},2) = \frac{D\left(\e^{\frac{2\pi\ii}{3}}\right)}{\pi}= \frac{\mathrm{Vol}(\mathrm{Fig}\, 8)}{2 \pi}\]

Here, $\mathrm{Vol} (\mathrm{Fig}\,8)$ is the hyperbolic volume of the Figure-8 knot complement.
The first equality was computed by Smyth \cite{Smyth:MahlerMeasure} and the last one was observed by Boyd \cite{Boyd:Millenium}.

It follows that for $\lambda$ sufficiently small we have:
\begin{eqnarray*}
\frac{\mathrm{Vol}(\mathrm{Fig}\, 8)}{\pi}&=& 2m(1+x+y)=m(Q_3)=m\left(\frac 1 {\lambda} (1-(1-\lambda Q_3))\right)\\
&=& -\log \lambda + m(1-(1-\lambda Q_3))\\
&=& -\log \lambda - \sum_{n=1}^{\infty} \frac {c_n}{n}\\
&=&-\log \lambda - \sum_{n=1}^{\infty} \frac 1 n \sum_{k=0}^{n} \binom n k (- \lambda)^k a_k\\
&=&-\log \lambda - \sum_{n=1}^{\infty} \frac 1 n \sum_{k=0}^{n} \sum_{j=0}^k \binom n k  {\binom {k} {j}}^2 \binom {2j}{j} (- \lambda)^k
\end{eqnarray*}
where $c_n:=[(1-\lambda Q_3)^n]_0$ and $a_n=[Q_3^n]_0$.

This type of expression for the volume of the Figure-8 knot should be compared to L\"uck's formula \cite{Lueck,Dasbach-Lalin} (Theorem 3). In that formula, the volume of a knot complement is expressed in terms of a similar formula, but the Mahler measure is computed over an element in the group ring of the fundamental group of the knot which is non-abelian. The coefficients in L\"uck's formula are notably hard to compute in practice.

\subsection{Predicting a recurrence relation from another recurrence relation}

Recall that we saw that the circuits in the triangular lattice and the even ones in the honeycomb lattice are related by
\[ b_n=\sum_{j=0}^n \binom{n}{j} (-3)^{n-j} a_j.\]
More generally, in the previous example, we have that
\[c_n=\sum_{j=0}^n \binom{n}{j} (-\lambda)^j a_j.\]
Given a recurrence relation for the $a_j$ we would like to obtain a recurrence relation for the $c_j$. The first observation is that we can assume that $-\lambda=1$, since it is easy to find a recurrence relation for $(-\lambda)^j a_j$ from the one for $a_j$. In other words, we can assume
\[c_n = \sum_{j=0}^n \binom{n}{j} a_j.\]
Consider one more time the generating function
\[u(t) = \sum_{n=0}^\infty a_n t^n.\]
For the $c_n$ we have
\[v(t) = \sum_{n=0}^\infty c_n t^n =\sum_{n=0}^\infty \sum_{j=0}^n \binom{n}{j}a_j  t^n.\]
If $|t|$ is small enough, we can invert the order in the sum,
\[v(t) =\sum_{j=0}^\infty a_j \sum_{n=0}^\infty \binom{n}{j}  t^n.\]
Now observe that
\begin{eqnarray*}
\sum_{n=0}^\infty \binom{n}{j}  t^n & =&  \frac{t^j}{j!} \sum_{n=0}^\infty n\cdots (n-j+1) t^{n-j} =\frac{t^j}{j!} \left(\frac{1}{1-t}\right)^{(j)} \\
&=& \frac{t^j}{j!} \frac{ j!}{(1-t)^{j+1}} = \frac{ t^j}{(1-t)^{j+1}}. 
\end{eqnarray*}
Putting everything together,
\[v(t) = \frac{1}{1-t}\sum_{j=0}^\infty a_j  \left(\frac{t}{1-t}\right)^j.\]
Now, a recurrence relation for $a_n$ is equivalent to a differential equation for $u(t)$, which translates into a differential equation for $v(t)$ and from there one obtains a recurrence for $c_n$.

As an example let us compute the recurrence in the case of $c_n$ with 
$\lambda = \frac{1}{6}$.

First, we set $b_n = \frac{a_n}{(-6)^n}$. We have
\[n^2a_n-(10n^2-10n+3)a_{n-1}+9(n-1)^2a_{n-2}=0.\]
Therefore,
\begin{equation}\label{dif}
12 n^2b_n+2(10n^2-10n+3)b_{n-1}+3(n-1)^2b_{n-2}=0,
\end{equation}
and
\[ c_n = \sum_{j=0}^n \binom{n}{j} b_j.\]
Let
\[ u(t)=\sum_{n=0}^\infty b_n t^n, \qquad v(t) = \sum_{n=0}^\infty c_n t^n 
= \frac{1}{1-t}\sum_{n=0}^\infty b_n \left( \frac{t}{1-t}\right)^n .\]
Set $s=\frac{t}{1-t}$, then
\[ (1-t)v(t) = u(s).\]

Notice that recurrence (\ref{dif}) translates into the differential equation
\[12t(tu''+u')+2t(10t^2u''+20t u'+u) +3 t^2( t^2u''+3tu'+u)=0.\]

Hence,
\[t^2(t+6)(3t+2)u''+t(9t^2+40t+12)u'+t(3t+2)u =0.\]

We have
\[ u'(s) =(1-t)^2((1-t)v'(t)-v(t)).\]

\[ u''(s) =(1-t)^3((1-t)^2v''(t)-4(1-t)v'(t)+2v(t)).\]

Replacing in
\[s^2(s+6)(3s+2)u''+s(9s^2+40s+12)u'+s(3s+2)u =0,\]
we obtain
\[t^2(t+2)(5t-6)((1-t)^2v''(t)-4(1-t)v'(t)+2v(t))\]
\[+t(19t^2 - 16t -12)((1-t)v'(t)-v(t))+t(t+2)v(t) =0.\]
 
Thus 
\[(5t^4 - t^3 - 16t^2 + 12t)v''(t)+ (20t^3 - 3t^2 - 32t + 12)v'(t)\]
\[+(10t^2 - t - 10)v(t)=0.\]

Finally, the differential equation translates into the recurrence
 \[12n^2c_n -2(8n^2-8 n+5)c_{n-1} - (n-1)^2 c_{n-2}+5(n-1)(n-2)c_{n-3}=0,\]
with initial terms $c_0=1, c_1=\frac{1}{2}, c_2=\frac{5}{12}$.

\subsection{Free groups}
If we reconsider the cases of $P_l= x_1+x_1^{-1}+ \dots +x_l+x_l^{-1}$ and $Q_l=(1+x_1+\dots +x_{l-1})(1+x_1^{-1}+ \dots +x_{l-1}^{-1})$ in the context of free variables, we have observed in  \cite{Dasbach-Lalin} that $a_n^{(P_l)}$ (respectively $a_n^{(Q_l)}$) counts the number of circuits of length $n$ (resp. $2n$) in a $2l$ (resp. $l$)-regular trees respectively. The generating function for circuits in a $d$-regular tree was computed by Bartholdi \cite{Bartholdi:Counting},
\[g_d(\lambda)=\frac{2(d-1)}{d-2+d\sqrt{1-4(d-1)\lambda^2}}.\]

For example, for $x+x^{-1}+y+y^{-1}$, one gets
\[u(\lambda) = \frac{3}{1+2\sqrt{1-12\lambda^2}}.\]  
If we let $v(\mu) = u(\lambda)$ with $\mu = \lambda^2$, then,
\[ (12\mu-1)(16\mu-1) v'' +2(240\mu - 19)v'+96v=0,\]
and $a_0=1$, $a_2=4$, 
\[na_{2n} -2(14n-9)a_{2n-2} +96(2n-3)a_{2n-4}=0.\] 

For the case of $(1+x+y)(1+x^{-1}+y^{-1})$, we have
\[u(\lambda) = \frac {4}{1+3\sqrt{1-8\lambda}}.\]  
In particular,
\[ (8\lambda-1)(9\lambda-1) u'' +2(90\lambda - 11)u'+36u=0,\]
and $a_0=1$, $a_1=3$,
\[na_{n} -(17n-12)a_{n-1} +36(2n-3)a_{n-2}=0.\] 

\subsection{A non-abelian, non-free group example}
Let us consider again the polynomial $P=x+x^{-1}+y+y^{-1}$ but this time with respect to  the group 
\[\Gamma=\left< x, y \, | \, x^2y=yx^2, y^2x=xy^2\right>.\]

\begin{figure}
\includegraphics[width=\textwidth,keepaspectratio]{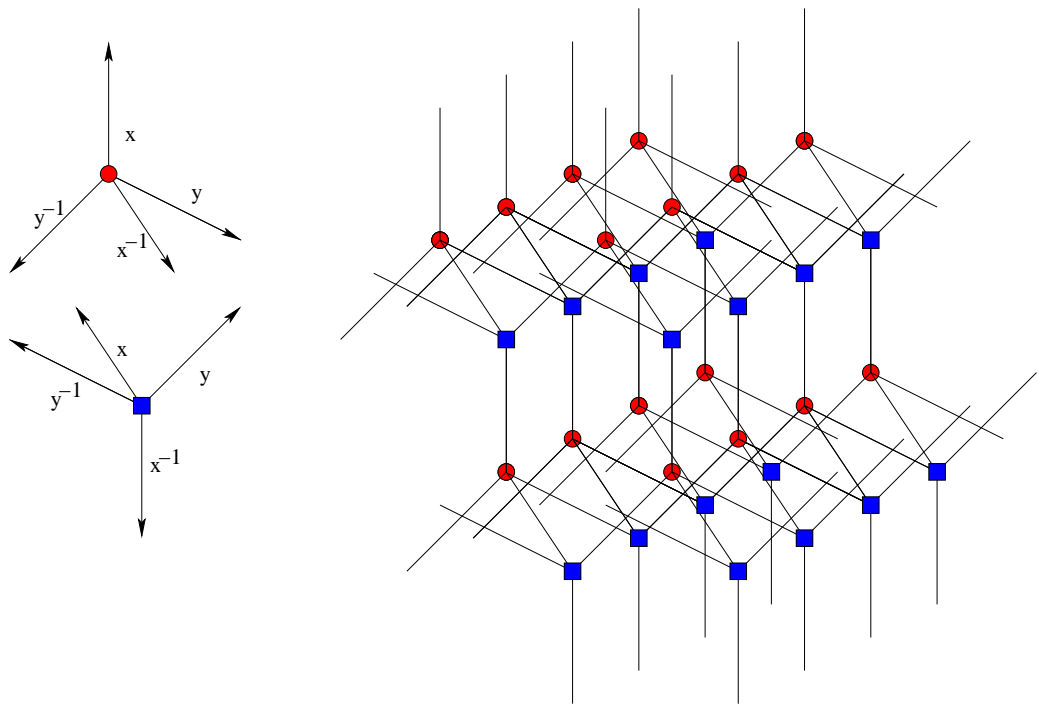}
\caption{The diamond lattice}
\label{Fig:diamonds}
\end{figure}

The $a_{2n}$ correspond to counting circuits of length $2n$ in the diamond lattice (Figure \ref{Fig:diamonds}). 

The vertices in this lattice can be divided into two 
different groups according to how the edges are oriented around the 
vertex. We indicate these two groups by rounded vertices and by a square vertices  in the picture. The 
models in the left show how to interpret a random walk that is leaving a 
vertex (rounded or square).

Notice that the edges can be divided into four families according to 
their direction. The minimum cycles in the lattice are given by hexagons 
and they correspond to the minimal relations in the group. Typically, a 
hexagon is formed by three pairs of parallel edges. Hence, there are four 
kinds of hexagons, according to which pair of parallel edges we choose to 
exclude. Now a simple inspection of the four cases of hexagons reveals 
that they all stand for either the relation $x^2y=yx^2$ or $y^2x=xy^2$ and 
those are the two generating relations.

The counting of the circuits in the diamond lattice appears in Domb \cite{Domb}. However, it is stated that the $a_n$ are the constant coefficients of powers of the polynomial
\[Q=(x+x^{-1}+z(y+y^{-1}))(x+x^{-1}+z^{-1}(y+y^{-1})),\]
with respect to the base group $\Z \times \Z \times \Z$. To see that both polynomials yield the same constant coefficients, one considers
\[ P^2= 4+x^2+x^{-2}+y^2+y^{-2}\]
\[+xy+xy^{-1}+x^{-1}y+x^{-1}y^{-1} +y^{-1}x^{-1}+ yx^{-1}+y^{-1}x+yx.\]
In $\Gamma$, even powers of $x$ and even powers of $y$ commute with each other and with the monomials $x^ay^b$ and $y^cx^d$, where each exponent is $1$ or $-1$. Also
\[  x^ay^by^cx^d= x^{a+d}y^{b+d},\]
because of the parity of the exponents.

On the other hand,
\[Q= 4+x^2+x^{-2}+y^2+y^{-2}\]
\[+z^{-1}xy+z^{-1}xy^{-1}+z^{-1}x^{-1}y+z^{-1}x^{-1}y^{-1} +zx^{-1}y^{-1}+ zx^{-1}y+zxy^{-1}+zxy.\]
Now we identify the monomials $x^ay^b$ (respectively $y^cx^d$) of $P^2$ with the monomials $z^{-1}x^ay^b$ (resp. $zx^dy^c$) of $Q$. It is an easy (and long) exercise to verify that corresponding monomials behave analogously in both cases.

With this interpretation, it is not hard to find a binomial formula for $a_n$ \cite{Domb}:
\[ a_{2n} = \sum_{k=0}^n \binom{n}{k}^2 \binom{2k}{k}\binom{2n-2k}{n-k}.\] 

The recurrence is given by $a_0=1$, $a_2=4$, and
\[ n^3a_{2n}-2(2n-1)(5n^2-5n+2)a_{2n-2}+64(n-1)^3a_{2n-4}=0.\]

Rogers \cite{Rogers} works with a third polynomial that yields the same Mahler measure:
\[R=4+\left(x+x^{-1}\right)\left(y+y^{-1}\right)+\left(y+y^{-1}\right)\left(z+z^{-1}\right)+\left(z+z^{-1}\right)\left(x+x^{-1}\right).\]  

An interesting fact is that the polynomial
\[S= \left(x+x^{-1}\right)\left(y+y^{-1}\right)+\left(y+y^{-1}\right)\left(z+z^{-1}\right)+\left(z+z^{-1}\right)\left(x+x^{-1}\right)\]  
corresponds to counting closed circuits in the face-centered cubic lattice. A closed form is given in \cite{Domb}:
\[b_n = \sum_{k_1 \equiv k_2 \equiv n \, \mathrm{mod 2}} \binom{n}{k_1+k_2} \binom{k_1+k_2}{k_1} \binom{k_1+k_2}{\frac{k_1+k_2}{2}} \binom{n-k_1}{\frac{n-k_1}{2}} \binom{n-k_2}{\frac{n-k_2}{2}}. \] 
We obtain
\[b_n=\sum_{j=0}^n \binom{n}{j} (-4)^{n-j} a_{2j}.\]

We compute the recurrence for $b_n$. First consider $a_n'=\frac{a_{2n}}{(-4)^n}$ and $b_n'=\frac{b_n}{(-4)^n}$. Then the $a_n'$ satisfy
\[2n^3a'_{n}+(2n-1)(5n^2-5n+2)a'_{n-1}+8(n-1)^3a'_{n-2}=0.\]
As before, we write
\[u(t) = \sum_{n=0}^\infty a'_n t^n \qquad v(t)= \sum_{n=0}^\infty b'_n t^n,\]
and $s=\frac{t}{1-t}$.

The differential equation for $u(s)$ is given by
\[(8s^5+10s^4+2s^3)u'''+(48s^4+45s^3+6s^2)u''+(56s^3+34s^2+2s)u'+(8s^2+2s)u=0.\]

In terms of $t$,
\[(6t^5 - 10t^4 + 2t^3 + 2t^2)v'''(t) +( 45t^4 - 60t^3 + 9t^2 + 6t)v''(t) \]
\[+(72t^3-72t^2+6t+2)v'(t)+(18t^2-12t)v(t)=0,\]
which translates into the recursion
\[2n^3b'_n+n(n-1)(2n-1)b'_{n-1}\]
\[-2(n-1)(5n^2-10n+6)b'_{n-2}+3(n-1)(n-2)(2n-3)b'_{n-3}=0.\]

Finally,
\[n^3b_n-2n(n-1)(2n-1)b_{n-1}\]
\[-16(n-1)(5n^2-10n+6)b_{n-2}-96(n-1)(n-2)(2n-3)b_{n-3}=0,\]
with $b_0=1, b_1=0, b_2=12$.

\section{Recurrences and hypergeometric functions}

A possible way to compute these recurrences is to use the algorithms 
in \cite{A=B}.

More explicitly, a result of Rogers \cite{Rogers} (based on techniques of Rodriguez-Villegas \cite{RV:ModularMahler}) relates the power series 
corresponding to the cubic lattice (easily related to the honeycomb as we 
have already noticed) and the diamond lattice to hypergeometric functions:

\begin{thm} (3.1 Rogers \cite{Rogers}) For $\lambda$ sufficiently small,

\begin{equation}
_3F_2\left(\frac{1}{3},\frac{1}{2},\frac{2}{3};1,1;-\frac{108 
\lambda}{(1-16\lambda)^3}\right) = (1-16 \lambda) \sum_{n=0}^ \infty 
\sum_{k=0}^n \binom{n}{k}^2 \binom{2k}{k} \binom{2n-2k}{n-k}\lambda^n
  \end{equation}

\begin{equation}
  _3F_2\left(\frac{1}{4},\frac{1}{2},\frac{3}{4};1,1;\frac{256 
\lambda}{9(1+3\lambda)^4}\right) = \frac{1+3\lambda}{1+\lambda} 
\sum_{n=0}^ \infty \binom{2n}{n} \sum_{k=0}^n \binom{n}{k}^2 
\binom{2k}{k}\left(\frac{\lambda}{9(1+\lambda)^2}\right)^n
  \end{equation}

\end{thm}

Here
\[_3F_2(a_1,a_2,a_3;b_1,b_2;\lambda)= \sum_{n=0}^\infty 
\frac{(a_1)_n(a_2)_n(a_3)_n}{(b_1)_n(b_2)_n} \frac{\lambda^n}{n!}\]
is a generalized hypergeometric series.

It satisfies the differential equation
\begin{eqnarray*}
\left((\vartheta_\lambda +b_1-1) (\vartheta_\lambda +b_2-1) - \lambda 
(\vartheta_\lambda +a_1)(\vartheta_\lambda +a_2)(\vartheta_\lambda 
+a_3)\right)  &&\\
{_3F_2}(a_1,a_2,a_3;b_1,b_2;\lambda) &=& 0,
\end{eqnarray*}
where $\vartheta_\lambda$ is the differential operator $\lambda \frac{d}{d 
\lambda}$.

It is possible then, to combine the differential equation for the 
generalized hypergeometric function and the formula for $u(\lambda)$ in 
order to obtain the recurrence. In fact, the details for the diamond 
lattice can be found in Section 4 of \cite{ChanChanLiu}.

{\bf Acknowledgements}: The authors would like to thank Neal Stoltzfus, Fernando Rodriguez-Villegas, and Pablo Bianucci for helpful discussions. ML expresses her gratitude to the Department of Mathematics at Louisiana State University for its hospitality.


 \small\parskip=-0.1cm

\renewcommand{\thefootnote}{}
\footnote{\null\hskip-0.65cm\null First author is supported by NSF-DMS-0456275 (FRG)}
\footnote{\null\hskip-0.65cm\null Second author is supported by University of Alberta Fac. Sci. startup grant N031000610}

\begin{flushright}
\addressa
\\
\emaila
\\[10pt]
\addressb
\\
\emailb
\end{flushright}

\label{lastpageR}
\end{document}